COOKED SUMMARY

This cooked summary explains the main ideas of the paper "**Max-Min SNR Signal Energy based Spectrum Sensing Algorithms for Cognitive Radio Networks with Noise Variance Uncertainty**" by referring the equation numbers in the paper. This cooked summary is very helpful for researchers in the field and who would like to know the novelty and key contribution of the paper without reading the whole text.

*Objective of the paper*

Given $N$ samples, the objective is to decide between $H_0$ (the $N$ samples contain noise only) and $H_1$ (these $N$ samples contain a transmitted signal + noise).

*Methodology*

- Express the baseband transmitted signal $x(t)$ as in Equation (1) of the paper. To do this we ASSUME that the transmitter pulse shaping filter is known to the cognitive receiver.
- Introduce a linear combination scalars $\{\alpha_i\}_{i=1}^L$ and DEFINE a new sample $\{\tilde{y}[n]\}_{n=1}^N$ as in Equation (3) (novel part).
- Design $\{\alpha_i\}_{i=1}^L$ such that we can get two different signals from $\{\tilde{y}[n]\}_{n=1}^N$ such that the SNR of the first signal is different from the SNR of the second signal. This is possible by performing the following tasks:
  - Solve the optimization problem (6), substitute this optimal solution $\{\alpha_i\}_{i=1}^L$ in Equation (3) and set the resulting samples as $\{e[n]\}_{n=1}^N$ (i.e., see Equation (14)) (novel part).
  - Solve the optimization problem (7), substitute this optimal solution $\{\alpha_i\}_{i=1}^L$ in Equation (3) and set the resulting samples as $\{z[n]\}_{n=1}^N$ (i.e., see Equation (14)) (novel part).
  - Now it is clear that the SNR of $z[n]$ is greater than (and is equal to) the SNR of $e[n]$ under $H_1$ (and $H_0$) hypothesis, respectively.
- Due to this mathematical outcome, we propose Equation (20) as our test statistics (novel part).
- As we can see, Equation (20) will be closer to 0 and much greater than 0 under $H_0$ and $H_1$ hypothesis, respectively.
- The $P_f$ and $P_d$ of this new test statistics is derived in (22) and (23).
- All the explanation and analytical equations after Equation (23) are to improve the test statistics (20) by taking into account the effect of synchronization between the transmitter and cognitive receiver, adjacent channel interference, unknown pulse shaping filter and so on, which are very important for practical cognitive radio.

# Max-Min SNR Signal Energy based Spectrum Sensing Algorithms for Cognitive Radio Networks with Noise Variance Uncertainty

Tadilo Endeshaw Bogale, *Student Member, IEEE* and Luc Vandendorpe, *Fellow, IEEE*

*Abstract*—This paper proposes novel spectrum sensing algorithms for cognitive radio networks. By assuming known transmitter pulse shaping filter, synchronous and asynchronous receiver scenarios have been considered. For each of these scenarios, the proposed algorithm is explained as follows: First, by introducing a combiner vector, an over-sampled signal of total duration equal to the symbol period is combined linearly. Second, for this combined signal, the Signal-to-Noise ratio (SNR) maximization and minimization problems are formulated as Rayleigh quotient optimization problems. Third, by using the solutions of these problems, the ratio of the signal energy corresponding to the maximum and minimum SNRs are proposed as a test statistics. For this test statistics, analytical probability of false alarm ($P_f$) and detection ($P_d$) expressions are derived for additive white Gaussian noise (AWGN) channel. The proposed algorithms are robust against noise variance uncertainty. The generalization of the proposed algorithms for unknown transmitter pulse shaping filter has also been discussed. Simulation results demonstrate that the proposed algorithms achieve better $P_d$ than that of the Eigenvalue decomposition and energy detection algorithms in AWGN and Rayleigh fading channels with noise variance uncertainty. The proposed algorithms also guarantee the desired $P_f(P_d)$ in the presence of adjacent channel interference signals.

## I. INTRODUCTION

The current wireless communication networks adopt fixed spectrum access strategy. The Federal Communications Commission have found that this fixed spectrum access strategy utilizes the available frequency bands inefficiently [1], [2]. A promising approach of addressing this problem is to deploy a cognitive radio (CR) network. One of the key characteristics of a CR network is its ability to discern the nature of the surrounding radio environment. This is performed by the spectrum sensing (signal detection) part of a CR network.

The most common spectrum sensing algorithms for CR networks are matched filter, energy and cyclostationary based algorithms. If the characteristics of the primary user such as modulation scheme, pulse shaping filter and packet format are known perfectly, matched filter is the optimal signal detection algorithm as it maximizes the received Signal-to-Noise Ratio

The authors would like to thank SES for the financial support of this work, the french community of Belgium for funding the ARC SCOOP and BELSPO for funding the IAP BESTCOM project. Part of this work has been published in the 8th International Conference on Cognitive Radio Oriented Wireless Networks (CROWNCOM), Washington DC, Jul. 8 - 10, 2013. Tadilo E. Bogale and Luc Vandendorpe are with the ICTEAM Institute, Université catholique de Louvain, Place du Levant 2, 1348 - Louvain La Neuve, Belgium. Email: {tadilo.bogale, luc.vandendorpe}@uclouvain.be, Phone: +3210478071, Fax: +3210472089.

(SNR). This algorithm has two major drawbacks: The first drawback is it needs dedicated receiver to detect each signal characteristics of a primary user [3]. The second drawback is it requires perfect synchronization between the transmitter and receiver which is impossible to achieve. This is due to the fact that, in general, the primary and secondary networks are administered by different operators. Energy detector does not need any information about the primary user and it is simple to implement. However, energy detector is very sensitive to noise variance uncertainty, and there is an SNR wall below which this detector can not guarantee a certain detection performance [3]–[5]. Cyclostationary based detection algorithm is robust against noise variance uncertainty and it can reject the effect of adjacent channel interference. However, the computational complexity of this detection algorithm is high, and large number of samples are required to exploit the cyclostationarity behavior of the received signal [5], [6]. On the other hand, this algorithm is not robust against cyclic frequency offset which can occur due to clock and timing mismatch between the transmitter and receiver [7]. In [8], Eigenvalue decomposition (EVD)-based spectrum sensing algorithm has been proposed. This algorithm is robust against noise variance uncertainty but its computational complexity is high. Furthermore, for single antenna receiver, this algorithm is sensitive to adjacent channel interference signal, and for multi-antenna receiver, this algorithm requires a channel covariance matrix different from a scaled identity [9].

This paper proposes novel spectrum sensing algorithms for cognitive radio networks. It is well known that a digital communication signal is constructed by passing an over-sampled signal through a transmitter pulse shaping filter. In most primary networks, for a given frequency band, as this pulse shaping filter is designed at the time of frequency planning stage and it is kept fixed, it is assumed to be known to the cognitive receiver. We consider synchronous and asynchronous receiver scenarios. For each of these scenarios, the proposed detection algorithm is explained as follows: First, by introducing a combiner vector, an over-sampled signal of total duration equal to the symbol period is combined linearly. Second, for this combined signal, the SNR maximization and minimization problems are formulated as Rayleigh quotient optimization problems. Third, by using the solutions of these problems, the ratio of the signal energy corresponding to the maximum and minimum SNRs are proposed as a test statistics. For this test statistics, analytical probability of false alarm ($P_f$) and probability of detection ($P_d$) expressions are derived for additive white Gaussian noise (AWGN) channel. The generalization of the proposed algorithms for unknown transmitter pulse shaping filter scenarios has also been discussed. It is shown that these detection algorithms (i.e., synchronous and asynchronous receiver scenarios) are robust against noise variance uncertainty. Moreover, under noise variance uncertainty, simulation results demonstrate that the proposed detection algorithms achieve better detection performance compared to that of the EVD-based and energy detection algorithms in AWGN and Rayleigh fading channels. The proposed detection algorithms also guarantee the desired $P_f(P_d)$ in the presence of low (moderate) adjacent channel interference (ACI) signals.

The remaining part of this paper is organized as follows: Section II discusses the hypothesis test problem. Section III presents the proposed spectrum sensing algorithms for a narrow band signal with synchronous and asynchronous receiver scenarios. In Section IV, the extension of the proposed spectrum sensing algorithms for a wide band signal is discussed. In Section V, computer simulations are used to compare the performance of the proposed and existing spectrum sensing algorithms. Conclusions are presented in Section VI.

*Notations:* The following notations are used throughout this paper. Upper/lower case boldface letters denote matrices/column vectors. The $\mathbf{X}_{(n,n)}$, $\mathbf{X}_{(n,:)}$, $\mathbf{X}^T$ and $\mathbf{X}^H$ denote the $(n,n)$ element, $n$th row, transpose and conjugate transpose of $\mathbf{X}$, respectively. $\mathbf{I}_n(\mathbf{I})$ is an identity matrix of size $n \times n$ (appropriate size) and, $(.)^\star$, $\mathrm{E}\{.\}$, $|.|$ and $(.)^*$ denote optimal, expectation, absolute value and conjugate operators, respectively.

## II. PROBLEM FORMULATION

Assume that the transmitted symbols $s_n, \forall n$ are pulse shaped by a filter $g(t)$. After the digital to analog converter, the baseband transmitted signal is given by

$$x(t) = \sum_{k=-\infty}^{\infty} s_k g(t - kT_s) \quad (1)$$

where $T_s$ is the symbol period. We assume that $x(t)$ is narrow band signal[1]. In an AWGN channel, the baseband received signal is expressed as

$$r(t) = \int_{-\infty}^{\infty} f^*(\tau)(x(t-\tau) + w(t-\tau))d\tau$$

$$= \int_{-\infty}^{\infty} f^*(\tau)(\sum_{k=-\infty}^{\infty} s_k g(t - kT_s - \tau) + w(t-\tau))d\tau$$

$$= \sum_{k=-\infty}^{\infty} s_k h(t - kT_s) + \int_{-\infty}^{\infty} f^*(\tau)w(t-\tau)d\tau$$

where $f^*(t)$ is the receiver filter, $w(t)$ is the additive white Gaussian noise and $h(t) = \int_{-\infty}^{\infty} f^*(\tau)g(t-\tau)d\tau$. The objective of spectrum sensing is to decide between $H_0$ and $H_1$

[1]The extension of the proposed spectrum sensing algorithms for a wide band signal will be discussed later.

from $r(t)$, where

$$r(t) = \int_{-\infty}^{\infty} f^*(\tau)w(t-\tau)d\tau, \quad H_0 \quad (2)$$

$$= \sum_{k=-\infty}^{\infty} s_k h(t - kT_s) + \int_{-\infty}^{\infty} f^*(\tau)w(t-\tau)d\tau, \; H_1.$$

Without loss of generality, we assume that $r(t)$ is a zero mean signal. Note that when $r(t)$ has a nonzero mean, its mean can be removed before examined by the proposed spectrum sensing algorithms.

## III. PROPOSED SPECTRUM SENSING ALGORITHMS

We define the $n$th discrete signal $\{\tilde{y}[n]\}_{n=1}^N$ as follows:

$$\tilde{y}[n] \triangleq \sum_{i=0}^{L-1} \alpha_i r((n-1)T_s + t_i) \quad (3)$$

$$= \sum_{k=-\infty}^{\infty} s_k \sum_{i=0}^{L-1} \alpha_i h((n-1)T_s + t_i - kT_s) +$$

$$\sum_{i=0}^{L-1} \alpha_i \int_{-\infty}^{\infty} f^*(\tau) w((n-1)T_s + t_i - \tau)d\tau$$

where $\{t_i\}_{i=0}^{L-1}$ are chosen such that $t_L - t_0 = T_s$ and $\{\alpha_i\}_{i=0}^{L-1}$ are the introduced variables. By assuming that the signal and noise (i.e., $x(t)$ and $w(t)$) are independent, the power of $\tilde{y}[n]$ can be expressed as

$$\mathrm{E}\{|\tilde{y}[n]|^2\} =$$

$$\mathrm{E}\{|\sum_{k=-\infty}^{\infty} s_k \sum_{i=0}^{L-1} \alpha_i h((n-1)T_s + t_i - kT_s)|^2\} +$$

$$\mathrm{E}\{|\sum_{i=0}^{L-1} \alpha_i \int_{-\infty}^{\infty} f^*(\tau) w((n-1)T_s + t_i - \tau)d\tau|^2\}$$

$$= \sigma_s^2 \sum_{k=-\infty}^{\infty} \boldsymbol{\alpha}^H \mathbf{A}_{nk} \boldsymbol{\alpha} + \sigma_w^2 \boldsymbol{\alpha}^H \mathbf{B}_n \boldsymbol{\alpha}$$

$$= \sigma_s^2 \boldsymbol{\alpha}^H \mathbf{A}_n \boldsymbol{\alpha} + \sigma_w^2 \boldsymbol{\alpha}^H \mathbf{B}_n \boldsymbol{\alpha} \quad (4)$$

where $\sigma_s^2$ and $\sigma_w^2$ are the variances of the signal and noise, respectively, $\boldsymbol{\alpha} = [\alpha_0, \alpha_1, \cdots, \alpha_{L-1}]^T$, $\mathbf{A}_{nk} = \mathbf{a}_{nk}\mathbf{a}_{nk}^H$, $\mathbf{A}_n = \sum_{k=-\infty}^{\infty} \mathbf{A}_{nk}$ and $\mathbf{B}_n = \frac{1}{\sigma_w^2}\mathrm{E}\{\mathbf{b}_n\mathbf{b}_n^H\}$ with $\mathbf{a}_{nk} = [h((n-1)T_s + t_0 - kT_s), h((n-1)T_s + t_1 - kT_s), \cdots, h((n-1)T_s + t_{L-1} - kT_s)]^T$ and $\mathbf{b}_n = [\int_{-\infty}^{\infty} f^*(\tau)w((n-1)T_s + t_0 - \tau)d\tau, \int_{-\infty}^{\infty} f^*(\tau)w((n-1)T_s + t_1 - \tau)d\tau, \cdots, \int_{-\infty}^{\infty} f^*(\tau)w((n-1)T_s + t_{L-1} - \tau)d\tau]^T$.

The entries of $\mathbf{A}_n$ and $\mathbf{B}_n$ can further be expressed as $(\mathbf{A}_n)_{(i+1,j+1)} = \sum_{k=-\infty}^{\infty} h((n-1-k)T_s + t_i)h^*((n-1-k)T_s + t_j) = \sum_{k'=-\infty}^{\infty} h(k'T_s + t_i)h^*(k'T_s + t_j) \triangleq \mathbf{A}_{(i+1,j+1)}$ and $(\mathbf{B}_n)_{(i+1,j+1)} = \int_{-\infty}^{\infty} f^*(\tau)f(t_i - t_j + \tau)d\tau \triangleq \mathbf{B}_{(i+1,j+1)}$. It follows

$$\mathrm{E}\{|\tilde{y}[n]|^2\} = \sigma_s^2 \boldsymbol{\alpha}^H \mathbf{A}\boldsymbol{\alpha} + \sigma_w^2 \boldsymbol{\alpha}^H \mathbf{B}\boldsymbol{\alpha}. \quad (5)$$

For given $\mathbf{A}$ and $\mathbf{B}$, the SNR minimization and maximization problems of $\mathrm{E}\{|\tilde{y}[n]|^2\}$ can be expressed as

$$\min_{\boldsymbol{\alpha}_{min}} \frac{\sigma_s^2 \boldsymbol{\alpha}_{min}^H \mathbf{A} \boldsymbol{\alpha}_{min}}{\sigma_w^2 \boldsymbol{\alpha}_{min}^H \mathbf{B} \boldsymbol{\alpha}_{min}} \equiv \min_{\boldsymbol{\alpha}_{min}} \frac{\boldsymbol{\alpha}_{min}^H \mathbf{A} \boldsymbol{\alpha}_{min}}{\boldsymbol{\alpha}_{min}^H \mathbf{B} \boldsymbol{\alpha}_{min}} \quad (6)$$

$$\equiv \min_{\boldsymbol{\alpha}_{min}} \frac{\boldsymbol{\alpha}_{min}^H (\mathbf{A}+\mathbf{B}) \boldsymbol{\alpha}_{min}}{\boldsymbol{\alpha}_{min}^H \mathbf{B} \boldsymbol{\alpha}_{min}}$$

$$\max_{\boldsymbol{\alpha}_{max}} \frac{\boldsymbol{\alpha}_{max}^H (\mathbf{A}+\mathbf{B}) \boldsymbol{\alpha}_{max}}{\boldsymbol{\alpha}_{max}^H \mathbf{B} \boldsymbol{\alpha}_{max}}. \quad (7)$$

These optimization problems are Rayleigh quotient problems. Since $\mathbf{A}$ and $\mathbf{B}$ are positive semidefinite matrices, the Generalized eigenvalue solution approach can be applied to get the optimal solutions of these problems which is summarized as follows [10], [11]:

As $\mathbf{B}$ is a positive semidefinite matrix, applying eigenvalue decomposition gives us

$$\mathbf{B} = \mathbf{U} \begin{pmatrix} \boldsymbol{\Sigma} & \mathbf{0} \\ \mathbf{0} & \mathbf{0} \end{pmatrix} \mathbf{U}^H \triangleq \mathbf{U}\mathbf{D}\mathbf{D}\mathbf{U}^H \quad (8)$$

where $\boldsymbol{\Sigma}$ is a diagonal matrix containing nonzero eigenvalues of $\mathbf{B}$, $\mathbf{U}$ is a unitary matrix and

$$\mathbf{D} = \begin{pmatrix} \boldsymbol{\Sigma}^{\frac{1}{2}} & \mathbf{0} \\ \mathbf{0} & \mathbf{0} \end{pmatrix}. \quad (9)$$

The pseudoinverse of $\mathbf{B}$ is given by

$$\mathbf{B}^\dagger = \mathbf{U} \begin{pmatrix} \boldsymbol{\Sigma}^{-1} & \mathbf{0} \\ \mathbf{0} & \mathbf{0} \end{pmatrix} \mathbf{U}^H = \mathbf{U}\tilde{\mathbf{D}}\tilde{\mathbf{D}}\mathbf{U}^H \quad (10)$$

where

$$\tilde{\mathbf{D}} = \begin{pmatrix} \boldsymbol{\Sigma}^{-\frac{1}{2}} & \mathbf{0} \\ \mathbf{0} & \mathbf{0} \end{pmatrix}. \quad (11)$$

By employing (8) - (11), and defining $\tilde{\boldsymbol{\alpha}} \triangleq \mathbf{D}\mathbf{U}^H \boldsymbol{\alpha}_{min}$ for (6) and $\tilde{\tilde{\boldsymbol{\alpha}}} \triangleq \mathbf{D}\mathbf{U}^H \boldsymbol{\alpha}_{max}$ for (7), we can rewrite the problems (6) and (7) as

$$\min_{\tilde{\boldsymbol{\alpha}}} \frac{\tilde{\boldsymbol{\alpha}}^H \tilde{\mathbf{A}} \tilde{\boldsymbol{\alpha}}}{\tilde{\boldsymbol{\alpha}}^H \tilde{\boldsymbol{\alpha}}} \quad (12)$$

$$\max_{\tilde{\tilde{\boldsymbol{\alpha}}}} \frac{\tilde{\tilde{\boldsymbol{\alpha}}}^H \tilde{\mathbf{A}} \tilde{\tilde{\boldsymbol{\alpha}}}}{\tilde{\tilde{\boldsymbol{\alpha}}}^H \tilde{\tilde{\boldsymbol{\alpha}}}} \quad (13)$$

where $\tilde{\mathbf{A}} = (\mathbf{U}\tilde{\mathbf{D}})^H(\mathbf{A}+\mathbf{B})(\mathbf{U}\tilde{\mathbf{D}}) = [\mathbf{I}\ \mathbf{0};\mathbf{0}\ \mathbf{0}]+\tilde{\mathbf{D}}\mathbf{U}^H\mathbf{A}\mathbf{U}\tilde{\mathbf{D}}$. The optimal $\tilde{\boldsymbol{\alpha}}$ and $\tilde{\tilde{\boldsymbol{\alpha}}}$ of these problems are given by the eigenvectors corresponding to the minimum and maximum nonzero eigenvalues of $\tilde{\mathbf{A}}$, respectively. Since $\tilde{\mathbf{A}}$ is also a positive semidefinite matrix, its minimum and maximum nonzero eigenvalues are always positive. The optimal solutions of the original problems (6) and (7) are thus given by $\boldsymbol{\lambda} \triangleq \boldsymbol{\alpha}_{min}^\star = \mathbf{U}\tilde{\mathbf{D}}\tilde{\boldsymbol{\alpha}}^\star$ and $\boldsymbol{\tau} \triangleq \boldsymbol{\alpha}_{max}^\star = \mathbf{U}\tilde{\mathbf{D}}(\tilde{\tilde{\boldsymbol{\alpha}}})^\star$.

At optimality, the denominator terms of the above problems are equal to unity (or any other positive value). Thus, under $H_0$ hypothesis, the optimal values of (12) and (13) are the same and equal to unity. However, under $H_1$ hypothesis, the optimal value of (13) is higher than that of (12)[2]. Due to this fact, we propose the following test statistics:

$$\widehat{\widetilde{T}} = \frac{\sum_{n=1}^{N} |\tilde{y}[n]|^2_{\boldsymbol{\alpha}_{max}}}{\sum_{n=1}^{N} |\tilde{y}[n]|^2_{\boldsymbol{\alpha}_{min}}} \triangleq \frac{\sum_{n=1}^{N} |z[n]|^2}{\sum_{n=1}^{N} |e[n]|^2} \triangleq \frac{\widehat{M}_{a2z}}{\widehat{M}_{a2e}} \quad (14)$$

where

$$\widehat{M}_{a2z} = \frac{1}{N} \sum_{n=1}^{N} |z[n]|^2, \ \widehat{M}_{a2e} = \frac{1}{N} \sum_{n=1}^{N} |e[n]|^2.$$

The authors of [8] propose over-sampling along with pre-whitening method to apply the EVD-based detection algorithm for the single receiver antenna case. However, in practice, there is always a nonzero (with very small power) adjacent channel interference signal. And, as will be clear in the simulation section, the algorithm of [8] can not ensure a predefined $P_f$ when there is an adjacent channel interference signal. However, as we can see from (14), the proposed test statistics can guarantee a predefined $P_f(P_d)$ when the adjacent channel interference signal power is very small compared to that of the desired signal and noise powers (see also the simulation section).

For sufficiently large $N$ (which is the case in a CR), by applying central limit theorem, we can interpret $z[n]$ and $e[n]$ as filtered and down-sampled versions of $\{w[i]\}_{i=1}^{LN}$, where the filters are $\boldsymbol{\eta}^{R+L-1\times 1} \triangleq \sqrt{2(1+\gamma_{max})} \sum \mathrm{diag}(\boldsymbol{\Upsilon},k)$ and $\boldsymbol{\theta}^{R+L-1\times 1} \triangleq \sqrt{2(1+\gamma_{min})} \sum \mathrm{diag}(\boldsymbol{\Psi},k)$ for $z[n]$ and $e[n]$, respectively, $k = [-(L-1), -(L-2), \cdots, R-1]$ and $\gamma_{max}$ and $\gamma_{min}$ denote the SNRs obtained by solving the problems (6) and (7), respectively, with $w[i], \forall i$ are independent and identically distributed (i.i.d) zero mean circularly symmetric complex Gaussian (ZMCSCG) random variables all with unit variance[3], $\boldsymbol{\Upsilon} = \boldsymbol{\tau}^T \otimes \mathbf{f}$, $\boldsymbol{\Psi} = \boldsymbol{\lambda}^T \otimes \mathbf{f}$, $\otimes$ denotes a kronecker product, $\mathbf{f} = [f_0, f_1, \cdots, f_R]$ is the sampled version of the receiver filter $f(t)$ with sampling period $\frac{T_s}{L}$, $R$ is the filter length and $\sum \mathrm{diag}(\mathbf{X},k)$ denotes the sum of the $k$th ($k = 0, k > 0$ and $k < 0$, denote the main diagonal, above the main diagonal and below the main diagonal, respectively) diagonal elements of $\mathbf{X}$.

For better exposition, let us introduce a new variable $\widetilde{T}$

$$\widetilde{T} = \frac{\lim_{N\to\infty} \frac{1}{N} \sum_{n=1}^{N} |z[n]|^2}{\lim_{N\to\infty} \frac{1}{N} \sum_{n=1}^{N} |e[n]|^2} \triangleq \frac{M_{a2z}}{M_{a2e}}. \quad (15)$$

By defining $\sigma_z^2 \triangleq 1+\gamma_{max}$, $\sigma_e^2 \triangleq 1+\gamma_{min}$ and $\gamma_d \triangleq \gamma_{max} - \gamma_{min}$, $\widetilde{T}$ can be expressed as

$$\widetilde{T} = \frac{2\sigma_z^2}{2\sigma_e^2} = 1, \qquad \mathrm{H}_0$$
$$= \frac{2\sigma_z^2}{2\sigma_e^2} = 1 + \frac{\gamma_d}{1+\gamma_{min}}, \qquad \mathrm{H}_1. \quad (16)$$

From this equation, one can notice that our problem turns to examining whether $\widetilde{T} = 1$ or $\widetilde{T} > 1$ for sufficiently large $N$. To get the $P_d$ and $P_f$ of the proposed test statistics, we examine the following Theorem [12].

---

[2]Note that under $H_1$ hypothesis, the optimal values of (12) and (13) are equal if and only if $\mathbf{A} = \rho\mathbf{B}$, where $\rho$ is any real number, which will never happen in a practical scenario.

[3]This is due to the fact that the noise power does not have any effect on the test statistics under $H_0$ hypothesis, and the effect of the signal power is incorporated by the filters $\boldsymbol{\eta}$ and $\boldsymbol{\theta}$ under $H_1$ hypothesis.

*Theorem 1*: Given a real valued function $\widehat{\widetilde{T}} = \frac{\widehat{M}_{a2z}}{\widehat{M}_{a2e}}$, the asymptotic distribution of $\sqrt{N}(\widehat{\widetilde{T}} - \widetilde{T})$ is given by

$$\sqrt{N}(\widehat{\widetilde{T}} - \widetilde{T}) \sim \mathcal{N}(0, \tilde{\sigma}^2) \quad (17)$$

where $\tilde{\sigma}^2 = \mathbf{v}^T \mathbf{\Phi} \mathbf{v}$,

$$\mathbf{v} = \left[\frac{\partial \widehat{\widetilde{T}}}{\partial \widehat{M}_{a2z}}, \frac{\partial \widehat{\widetilde{T}}}{\partial \widehat{M}_{a2e}}\right]^T_{\widehat{M}_{a2z}=M_{a2z}, \widehat{M}_{a2e}=M_{a2e}}$$
$$= \left[\frac{1}{M_{a2e}}, -\frac{M_{a2z}}{M_{a2e}^2}\right]^T \quad (18)$$

and $\mathbf{\Phi}$ is the asymptotic covariance matrix of a multivariate random variable $\sqrt{N}([\widehat{M}_{a2z}, \widehat{M}_{a2e}]^T - [M_{a2z}, M_{a2e}]^T) \sim \mathcal{N}(\mathbf{0}, \mathbf{\Phi})$.

*Proof:* See *Theorem* 3. 3. A on page 122 of [12]. ∎

Substituting $\mathbf{\Phi}$ into (17) gives

$$\tilde{\sigma}^2 = \frac{M_{a2e}^2 \mathbf{\Phi}_{(1,1)} - 2M_{a2e}M_{a2z}\mathbf{\Phi}_{(1,2)} + M_{a2z}^2 \mathbf{\Phi}_{(2,2)}}{M_{a2e}^4}. \quad (19)$$

The coefficients of $\mathbf{\Phi}$ can be computed numerically (see Appendix A).

As $\widetilde{T} = 1$ under $H_0$ hypothesis, we modify the test statistics $\widehat{\widetilde{T}}$ to

$$T = \sqrt{N}(\widehat{\widetilde{T}} - 1). \quad (20)$$

The $P_f$ of this test statistics is expressed as

$$P_f(\lambda) = Pr\{T > \lambda | H_0\}. \quad (21)$$

Under $H_0$ hypothesis, as $T \sim \mathcal{N}(0, \tilde{\sigma}_{H0}^2)$, the $P_f$ is given by

$$P_f = \int_\lambda^\infty \frac{1}{\sqrt{2\pi\tilde{\sigma}_{H0}^2}} \exp^{-\frac{x^2}{2\tilde{\sigma}_{H0}^2}} dx$$
$$= Q\left(\frac{\lambda}{\tilde{\sigma}_{H0}}\right) \quad (22)$$

where $Q(.)$ is the Q-function which is defined as [13]

$$Q(\lambda) = \frac{1}{\sqrt{2\pi}} \int_\lambda^\infty \exp^{-\frac{x^2}{2}} dx$$

and $\tilde{\sigma}_{H0}^2$ is $\tilde{\sigma}^2$ of (19) under $H_0$ hypothesis.

Mathematically, $P_d(\lambda)$ is expressed as

$$P_d(\lambda) = Pr\{T > \lambda | H_1\} \quad (23)$$
$$= \int_\lambda^\infty \frac{1}{\sqrt{2\pi\tilde{\sigma}_{H1}^2}} \exp^{-\frac{(x-\mu)^2}{2\tilde{\sigma}_{H1}^2}} dx = Q\left(\frac{\lambda - \mu}{\tilde{\sigma}_{H1}}\right)$$

where $\mu = \sqrt{N}(\widetilde{T} - 1) = \sqrt{N}\frac{\gamma_d}{1+\gamma_{min}}$ and $\tilde{\sigma}_{H1}^2$ is $\tilde{\sigma}^2$ of (19) under $H_1$ hypothesis. From the above expression, we can understand that for given $\gamma_d > 0$ and $\lambda$, increasing $N$ increases $P_d$. This is due to the fact that $Q(.)$ is a decreasing function. Thus, the proposed detection algorithm is consistent (i.e., for any given $P_f > 0$ and SNR, as $N \to \infty, P_d \to 1$).

As can be seen from (6) and (7), for a given $g(t)$, the achievable maximum and minimum SNRs depend on the selection of $f(t)$, $L$ and $\{t_i\}_{i=0}^{L-1}$. For a given $g(t)$, getting the optimal $f(t)$, $L$ and $\{t_i\}_{i=0}^{L-1}$ ensuring the highest detection performance is an open research topic. In our simulation, we have observed better detection performance when we select $f(t) = g(t)$ (i.e., matched filter), $L \geq 8$ and $\{t_i = T_s(\frac{1}{2} + \frac{i}{L})\}_{i=0}^{L-1}$. For example, if $f(t)$ is square root raised cosine filter (SRRCF), the initial timing ($t_0$) will be as in Fig. 1.

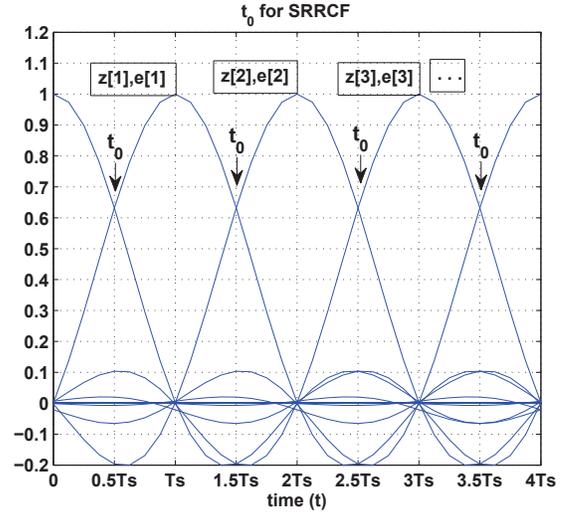

Fig. 1. Description of $t_0$ for SRRCF.

From (23) we can also notice that $P_d$ increases as $\frac{\gamma_d}{1+\gamma_{min}}$ increases. This is achieved when $\gamma_d > 0$. To get $\gamma_d = \gamma_{max} - \gamma_{min} > 0$, the ranks of $\mathbf{A}$ and $\mathbf{B}$ must be at least 2. When $f(t)$ is SRRCF with roll-off factors 0.2, 0.25 and 0.35, we have observed that $\frac{\gamma_d}{1+\gamma_{min}}$ increases as the roll-off factor increases, whereas, the ranks of $\mathbf{A}$ and $\mathbf{B}$ are the same (which is equal to 4) for these roll-off factors. From this explanation, we can understand that the quality of the proposed detector can not be determined just from the ranks of $\mathbf{A}$ and $\mathbf{B}$.

From Fig. 1, one can realize that to get the $P_d$ of (23), $t_0$ must be known perfectly. The exact $t_0$ is known when the receiver is synchronized perfectly with the transmitter. However, in general, since the transmitters and receivers are administered by different operators, perfect synchronization is not possible. Even in some scenario, the pulse shaping filter may not be known to the cognitive receiver. In the following, we generalize the aforementioned detector for known and unknown pulse shaping filters with asynchronous receiver scenarios.

### A. Asynchronous receiver and known transmitter pulse shaping filter scenario

In this subsection, we generalize the aforementioned detection algorithm for known pulse shaping with asynchronous receiver scenario. This algorithm is designed based on the estimation of $t_0$. As can be seen from (3), there are $L$ possible values of $t_0$. Thus, from the received signal $r(t)$, we can obtain $L$ possible values of $\widehat{\widetilde{T}}(T)$. Under $H_0$ hypothesis, all values of $\{T_i\}_{i=1}^L$ are almost the same, whereas, under $H_1$ hypothesis, the values of $\{T_i\}_{i=1}^L$ are not the same. And, the $t_0$ corresponding to $T_{max} = \max[T_1, T_2, \cdots, T_L]$ can be considered as the best estimate of the true $t_0$. Due to this

reason, we propose the following test statistics:

$$T_{max} = \max[T_1, T_2, \cdots, T_L]. \quad (24)$$

The cumulative distribution function (CDF) of $T_{max}$ is given by [14]

$$F_{T_{max}}(\tilde{\lambda}) = \int_{-\infty}^{\tilde{\lambda}} \int_{-\infty}^{\tilde{\lambda}} \cdots \int_{-\infty}^{\tilde{\lambda}} \frac{1}{(2\pi)^{L/2}|\Sigma|^{1/2}} \times \exp^{-\frac{1}{2}(\omega-\mu)^T\Sigma^{-1}(\omega-\mu)} d\omega_1 d\omega_2 \cdots d\omega_L \quad (25)$$

where $\mu$ and $\Sigma$ are the mean vector and covariance matrix of the random variables $[T_1, T_2, \cdots, T_L]$, respectively and $\omega = [\omega_1, \omega_2, \cdots, \omega_L]^T$. By defining $\{\mathbf{t}_k \triangleq [T_k, T_{k+1}, \cdots, T_L, T_1, \cdots, T_{k-1}] \sim (\mu_k, \Sigma_k)\}_{k=1}^L$, it can be easily seen that $\{\mu_k, \Sigma_k\}_{k=1}^L$ are not necessarily the same if the mean and variances of $\{T_i\}_{i=1}^L$ are not the same.

When the receiver is not synchronized with the transmitter, we can obtain $L$ possible values of $\mu$ and $\Sigma$ (i.e., $\{(\mu_k, \Sigma_k)\}_{k=1}^L$) all with equal probability of occurrence. Thus, the $P_f$ and $P_d$ of the test statistics $T_{max}$ can be expressed as

$$P_f = \frac{1}{L}\sum_{k=1}^{L} Pr\{T_{max} > \tilde{\lambda}|(H_0, \mathbf{t}_k)\}$$

$$= \frac{1}{L}\sum_{k=1}^{L}(1 - F_{T_{max}}|(H_0, \mathbf{t}_k)(\tilde{\lambda}))$$

$$= 1 - F_{T_{max}}|(H_0, \mathbf{t}_1)(\tilde{\lambda})$$

$$P_d = \frac{1}{L}\sum_{k=1}^{L} Pr\{T_{max} > \tilde{\lambda}|(H_1, \mathbf{t}_k)\}$$

$$= \frac{1}{L}\sum_{k=1}^{L}(1 - F_{T_{max}}|(H_1, \mathbf{t}_k)(\tilde{\lambda}))$$

$$= 1 - \frac{1}{L}\sum_{k=1}^{L} F_{T_{max}}|(H_1, \mathbf{t}_k)(\tilde{\lambda}) \quad (26)$$

where $F_{T_{max}}|(H_0, \mathbf{t}_k)(\tilde{\lambda})$ and $F_{T_{max}}|(H_1, \mathbf{t}_k)(\tilde{\lambda})$ are (25) under $H_0$ and $H_1$ hypothesis with the statistics of $\mathbf{t}_k$, respectively. The third equality follows from the fact that $\{\mathbf{t}_k \sim (\mu, \Sigma)\}_{k=1}^L$ under $H_0$ hypothesis. For a given $\mathbf{t}_k$, $\mu_k$ can be computed like $\mu$ of (23) and $\Sigma_k$ can be computed numerically (see Appendix B).

To the best of our knowledge, there is no any closed form solution for the integral (25). Due to this reason, this paper solves (25) using "mvncdf" matlab numerical package[4]. This package computes $F_{T_{max}}(\tilde{\lambda})$ for a given $\tilde{\lambda}$.

In practice, however, we are required to get $\tilde{\lambda}$ for a desired $P_f$ (i.e., a constant false alarm rate detector). In this paper, we apply a bisection search method to get $\tilde{\lambda}$ satisfying the $P_f$ of (26) [15]. To apply the bisection search, the lower and upper bounds of $\tilde{\lambda}$ are required which can be computed as follows:

Since $T_{max}$ is the supreme of all $\{T_i\}_{i=1}^L$, for fixed $P_f$, one can notice that $\tilde{\lambda} \geq \lambda = \tilde{\lambda}_{min}$, where $\lambda$ is as given in

[4]Note that this package is designed to compute the CDF of a multivariate Gaussian random variable.

(22). On the other hand, $\tilde{\lambda}$ of (26) becomes maximum when $\{T_i\}_{i=1}^L$ are independent. In such a case, the exact closed form expression of $F_{T_{max}}(\tilde{\lambda})$ is given as [14]

$$F_{T_{max}}(\tilde{\lambda}) = (F_T(\tilde{\lambda}))^L \quad (27)$$

where $F_T(\tilde{\lambda}) = 1 - Q\left(\frac{\tilde{\lambda}}{\tilde{\sigma}_{H0}}\right)$ is the CDF of $T$ (20). Thus, the maximum possible value of $\tilde{\lambda}$ is given by

$$P_f = 1 - F_{T_{max}}(\tilde{\lambda}_{max}) = 1 - (F_{T_1}(\tilde{\lambda}_{max}))^L$$

$$\Rightarrow Q\left(\frac{\lambda}{\tilde{\sigma}_{H0}}\right) = 1 - \left(1 - Q\left(\frac{\tilde{\lambda}_{max}}{\tilde{\sigma}_{H0}}\right)\right)^L$$

$$\Rightarrow \tilde{\lambda}_{max} = \tilde{\sigma}_{H0}Q^{-1}\left(1 - \left(1 - Q\left(\frac{\lambda}{\tilde{\sigma}_{H0}}\right)\right)^{1/L}\right) \quad (28)$$

where the third equality is due to (22).

The bisection search method for computing the exact $\tilde{\lambda}$ is summarized in **Algorithm I**.

**Algorithm I**
Initialization: Set $\tilde{\lambda}_{min}$ and $\tilde{\lambda}_{max}$ as in (22) and (28), respectively and $\varepsilon = 10^{-3}$.
**Repeat**:
Set $\bar{\lambda} = \frac{1}{2}(\tilde{\lambda}_{min} + \tilde{\lambda}_{max})$.
1) Compute $\bar{P}_f = 1 - F_{T_{max}}(\bar{\lambda})$ by employing "mvncdf" matlab package.
2) If $\bar{P}_f \leq P_f$, set $\tilde{\lambda}_{min} = \bar{\lambda}$ else set $\tilde{\lambda}_{max} = \bar{\lambda}$
**Until** $|\bar{P}_f - P_f| \leq \varepsilon$.
Set $\tilde{\lambda} = \bar{\lambda}$ as the threshold.

*B. Asynchronous receiver and unknown transmitter pulse shaping filter scenario*

In this subsection, the generalization of the proposed algorithm for the detection of a signal with unknown pulse shaping filter is discussed. For this scenario, there are two obvious questions: The first question is how can we select the receiver filter. The second question is how can we optimize $\alpha$ to achieve the best detection performance in asynchronous receiver scenario[5]. To address these questions, let us examine the detection of DVB-S2 signals. According to [16], a DVB-S2 signal employs a SRRCF with roll-off factor $0.2, 0.25$ or $0.35$.

*Selection of the receiver filter (***B***)*: We design **B** based on the smallest roll-off factor (i.e., a SRRCF with roll-off factor $0.2$). This is due to the fact that if we design **B** with a roll-off factor $> 0.2$, the examined band contains strong and unknown adjacent channel interference signal when the roll-off factor of the transmitter filter is smaller than that of the receiver filter. Consequently, a predefined $P_f$ cannot be ensured under $H_0$ hypothesis.

*Optimization of $\alpha$*: We optimize $\alpha$ by considering all possible pulse shaping filters and taking into account the probability of each pulse shaping filter. The optimal $\alpha$ can

[5]As the transmitted signal pulse shaping filter is not known, perfect synchronization between the transmitter and receiver can never be achieved.

be obtained by solving the following problems:

$$\min_{\boldsymbol{\alpha}_{minU}} \sum_{p=1}^{\tilde{P}} \theta_p \frac{\boldsymbol{\alpha}_{minU}^H \mathbf{A}_p \boldsymbol{\alpha}_{minU}}{\boldsymbol{\alpha}_{minU}^H \mathbf{B} \boldsymbol{\alpha}_{minU}} \quad (29)$$

$$\max_{\boldsymbol{\alpha}_{maxU}} \sum_{p=1}^{\tilde{P}} \theta_p \frac{\boldsymbol{\alpha}_{maxU}^H \mathbf{A}_p \boldsymbol{\alpha}_{maxU}}{\boldsymbol{\alpha}_{maxU}^H \mathbf{B} \boldsymbol{\alpha}_{maxU}} \quad (30)$$

where $\mathbf{A}_p$ is the matrix $\mathbf{A}$ of (5) corresponding to the $p$th pulse shaping filter, $\theta_p$ is the probability of the $p$th pulse shaping filter and $\tilde{P}$ is the number of possible pulse shaping filters[6]. These two problems can be examined exactly like those of (6) and (7).

With the optimal $\boldsymbol{\alpha}_{maxU}(\boldsymbol{\alpha}_{minU})$ of the above problem, the test statistics and $P_f(P_d)$ expressions can be expressed like that of Section III-A (i.e., the detection algorithm with known pulse shaping filter and asynchronous receiver scenario). The details are omitted for conciseness.

As the algorithm of this subsection employs the same receiver filter and $\boldsymbol{\alpha}$ for all transmitted signals, its implementation cost is lower than that of the algorithm of Section III-A. However, the $P_d$ of the algorithm of this subsection cannot be higher than that of the algorithm of Section III-A. This is due to the fact that when the roll-off factor of the transmitter pulse shaping filter is higher than 0.2, the receiver also employs the same roll-off factor (i.e., $> 0.2$). This increases the effective bandwidth of the received signal which subsequently increase the number of received independent samples and detection probability. Therefore, the tradeoff between the detection algorithm of this subsection and that of Section III-A is complexity versus performance.

As can be seen from (4), the entries of $\mathbf{B}$ can be obtained analytically from $f(t)$, whereas, the entries of $\mathbf{A}$ are obtained by infinite summation (i.e., $-\infty \leq k' \leq \infty$). However, in a practical filter, as the magnitudes of $f(t)(g(t))$ decrease as $|t|$ increases, the coefficients of $\mathbf{A}$ can be well approximated by employing finite summations (i.e., $-K' \leq k' \leq K'$), where $K'$ is a finite integer. We have provided a strategy for determining $\mathbf{A}$, $\mathbf{B}$ and $\boldsymbol{\alpha}$, and numerical results are presented for the scenario where $g(t)$ is a SRRCF with roll-off 0.2 and L=8 (see Appendix C).

For any transmitter pulse shaping filter $g(t)$, the proposed detectors are summarized in **Algorithm II**.

**Algorithm II**

**Synchronous receiver scenario**
**Initialization:** Set f(t)=g(t) (matched filtering) and $L$ and $R$ as required.
a) Search $t_0$ such that $\gamma_d$ is maximum. We would like to mention here that for the well known SRRCF pulse shaping filter, we have found almost constant $\gamma_d$ for $L \geq 8$. In our simulation, we choose $L = 8$ to reduce the computational complexity of the detector. However, for general pulse shaping filter, exhaustive search of $t_0$ can be applied for any $L$. This is due to the fact that the optimization problems are solved only once prior to the detection process.
b) With the above $t_0$, solve the optimization problems (6) and (7), and compute $\gamma_{min}$ and $\gamma_{max}$.
c) With the optimal $\boldsymbol{\alpha}$ of (6) and (7), compute $P_f$ using (22).
d) With the above $\gamma_{min}$, $\gamma_{max}$ and optimal $\boldsymbol{\alpha}$ of (6) and (7), compute $P_d$ using (23).

**Asynchronous receiver with known pulse shaping filter scenario**
**Initialization:** Set f(t)=g(t) (matched filtering) and $L$ and $R$ as required.
a) With the optimal $\boldsymbol{\alpha}$ of (6) and (7), compute $P_f$ and $P_d$ using (26).

**Asynchronous receiver with unknown pulse shaping filter scenario**
**Initialization:** Set f(t) as SRRCF with the smallest of all transmitted signal roll-off factors, and $L$ and $R$ as required.
a) With the optimal $\boldsymbol{\alpha}$ of (29) and (30), compute $P_f$ and $P_d$ using (26).

## IV. EXTENSION OF THE PROPOSED ALGORITHMS FOR DETECTION OF WIDE-BAND SIGNALS

In this section, we discuss the extension of the proposed algorithms for the detection of a wide band signal (i.e., a signal with frequency selective channel). It is known that any wide band signal is a superposition of two or more narrow band signals. Thus, a wide band signal ($x_w(t)$) can be expressed as

$$x_w(t) = \sum_{c=1}^{C} x_c(t) \quad (31)$$

where $x_c(t) = \sum_{k=-\infty}^{\infty} s_{ck} g(t - kT_s)$ with $s_{ck}, \forall k$ are the transmitted symbols of the $c$th narrow band signal and $C$ is the number of narrow band signals. For the detection of a wide band signal, we can modify the test statistics (20) to

$$T_w = \sum_{c=1}^{C} \sqrt{N}(\widehat{\tilde{T}}_{wc} - 1) \quad (32)$$

where

$$\widehat{\tilde{T}}_{wc} = \frac{\sum_{n=1}^{N} |\tilde{y}_c[n]|^2_{\boldsymbol{\alpha}_{max}}}{\sum_{n=1}^{N} |\tilde{y}_c[n]|^2_{\boldsymbol{\alpha}_{min}}} \quad (33)$$

with $\tilde{y}_c[n]$ as the $c$th narrow band signal $n$th sample. As we can see, the $P_d$ and $P_f$ of the test statistics (32) can be obtained exactly like the $P_d$ and $P_f$ expressions of Section III. The details are omitted for conciseness.

## V. SIMULATION RESULTS

In this section, we provide simulation results. All of the simulation results are obtained by averaging 10000 experiments. The SNR is defined as $SNR \triangleq \frac{\sigma_s^2}{\sigma_w^2}$ and $N = 2^{15}$. For better exposition, we assume that the transmitted signal is QPSK, and the channel between the transmitter and receiver is flat fading (i.e., narrow-band signal). In all of the figures, "Simu" and "Theo" represent simulation and theory, respectively, "Sync" denotes the proposed algorithm for synchronous receiver scenario, "Async with est" denotes the proposed algorithm for asynchronous receiver scenario which is designed based on the estimation of $t_0$ (i.e., the algorithm in Section III-A) and "Async w/o est" denotes the algorithm of [17] for asynchronous receiver scenario which is designed without estimating $t_0$. Unless otherwise stated explicitly, the transmitter and receiver employ a SRRCF with roll-off factor 0.2, $L=8$, $R=64L+1$ and $|f(t)|^2 = |g(t)|^2 = 1$.

---
[6] In practice $\tilde{P}$ is small.

## A. Verification of the $P_f$ expressions

In this subsection, we verify the theoretical $P_f$ expressions by computer simulation. We also examine the effect of adjacent channel interference signal on the $P_f$ expressions of the proposed algorithms and that of the EVD-based detection algorithm [8][7]. We employ an adjacent channel signal $A(t) = \kappa(I_1 \tilde{a}_1(t) + I_2 \tilde{a}_2(t))$ with $\sigma_{\tilde{a}_1}^2 = \sigma_{\tilde{a}_2}^2 = \sigma_w^2$, where $\tilde{a}_1(t) = \sin(\frac{2.4\pi}{T_s}t)a_1(t)$, $\tilde{a}_2(t) = \sin(\frac{4.4\pi}{T_s}t)a_2(t)$, $\kappa$, $a_1(t)(a_2(t))$ is a zero mean pulse shaped (with SRRCF and roll-off = 0.2) binary phase shift keying (BPSK) signal with symbol period $T_s$, and $I_1(I_2)$ is a discrete random variable which takes a value 0 or 1. With this adjacent channel signal, we get $r(t) = \int_{-\infty}^{\infty} f^\star(\tau)(w(t-\tau) + A(t-\tau))d\tau$ under $H_0$ hypothesis[8].

The $P_f$ expressions of the proposed detection algorithms and the EVD-based detection algorithm of [8] are plotted in Fig. 2. As can be seen from this figure, in the proposed algorithms, the theoretical $P_f$ expressions fit that of the simulation results for both $A(t) = 0$ and $A(t) \neq 0$, whereas, in the EVD-based algorithm of [8], the theoretical $P_f$ expression is deviated significantly from the simulation result in the practically relevant regions (i.e., the regions $0 \leq P_f \leq 0.1$) when $A(t) \neq 0$. From this discussion, we can understand that the proposed detection algorithms can maintain the required $P_f$ even in the presence of ACI signal.

## B. Verification of the theoretical $P_f$ versus $P_d$ expressions

In this subsection, we verify the theoretical $P_f$ versus $P_d$ expressions of the proposed detectors in AWGN channel by computer simulations. It is assumed that $A(t) = 0$, and the SNR is known perfectly[9] and it is set to $-14dB$. As can be seen from Fig. 3, the theoretical $P_f$ versus $P_d$ expressions fit that of the simulation in both synchronous and asynchronous receiver scenarios.

[7]As the maximum to minimum Eigenvalue (MME) detection algorithm of [8] gives superior performance, for the comparison, we employ the MME detection algorithm of [8]. For this algorithm, we use a smoothing factor of $4L$, and the Tracy-Widom distribution of order 2 ($TW_2$) values are taken from Table 3 of [18].

[8]The bandwidth of $x(t)$ is $\frac{1}{T_s}(1 + \text{roll} - \text{off})$Hz. This shows that most of the energy of $A(t)$ does not lie in the bandwidth of $x(t)$. And, the actual interference power (IP) is $IP = \text{E}|\int_{-\infty}^{\infty} f^\star(\tau)(A(t-\tau))d\tau|^2 \neq 0$ and indeed $IP < \sigma_{\tilde{a}_1}^2$. From this explanation, we can understand that $\sigma_{\tilde{a}_1}^2$ and $\sigma_{\tilde{a}_2}^2$ are not the true interference powers in the band of interest.

[9]Here $A(t) = 0$ and the true SNR (i.e., accurate signal and noise variances) are required just to get $P_d$ which depends on SNR.

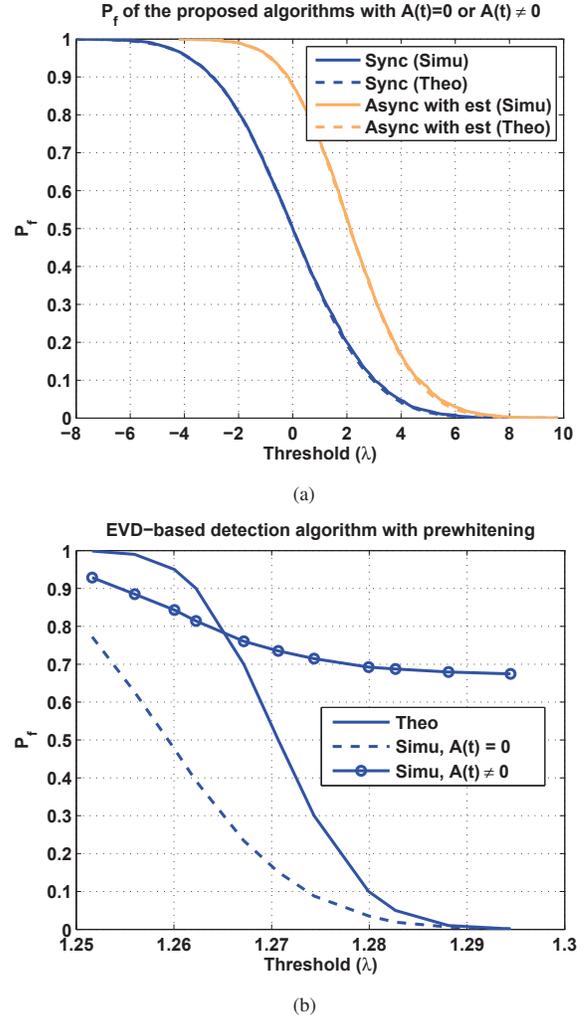

Fig. 2. Theoretical and simulated $P_f$ with and without ACI signal. (a) The proposed detection algorithms. (b) The EVD-based detection algorithm of [8].

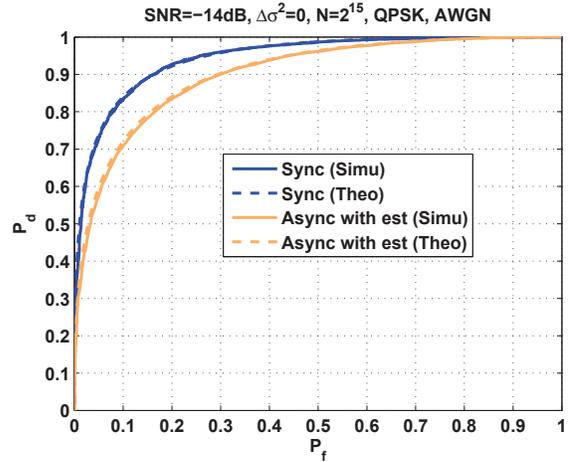

Fig. 3. Theoretical and simulated $P_f$ versus $P_d$ of the proposed detectors.

## C. Comparison of the proposed and existing detection algorithms under noise variance uncertainty

In this simulation, we compare the performance of the proposed detection algorithms (with and without adjacent channel interference signals) to those of [17], EVD-based and energy detection algorithms. For $A(t) \neq 0$, the adjacent channel interference signal of Section V-A with $\sigma_{\tilde{a}_1}^2 = \sigma_{\tilde{a}_2}^2 = \sigma_s^2$ is employed. The comparison is performed for AWGN and Rayleigh fading channels under noise variance uncertainty. According to [4], in an uncertain noise variance scenario, the true noise variance can be modeled as a bounded interval of $[\frac{1}{\epsilon}\sigma_w^2, \epsilon\sigma_w^2]$ for some $\epsilon = 10^{\Delta\sigma^2/10} > 1$, where the uncertainty $\Delta\sigma^2$ is expressed in dB. We assume that this bound follows a uniform distribution, i.e., $\mathcal{U}[\frac{1}{\epsilon}\sigma_w^2, \epsilon\sigma_w^2]$. The noise variance is the same for one experiment (since it has a short duration) and follows a uniform distribution during several experiments. Moreover, in a Rayleigh fading channel, the channel gain is the same for one experiment and follows a Rayleigh distribution during several experiments. The comparisons are performed for different SNRs by setting $P_f \leq 0.1$. Fig. 4 shows the performance of the proposed detection algorithms and those of the algorithm in [17], EVD-based and energy detection algorithms.

From Fig. 4, we observe that the proposed detection algorithms and that of [17] can maintain the desired $P_d$ in the presence of ACI signal and the best performance is achieved when the receiver is synchronized with the transmitter. Furthermore, the proposed algorithms (Sync and Async with est) achieve better detection performance compared to those of the algorithm in [17] (Async w/o est), EVD-based and energy detection algorithms.

We would like to mention here that when $A(t) \neq 0$, the exact thresholds for achieving $P_f \leq 0.1$ are not known for EVD-based and energy detection algorithms. For this reason, the $P_d$ curves of Fig. 4 employ $A(t) = 0$ only for EVD-based and energy detection algorithms.

## D. The effect of ACI on the $P_f$ and $P_d$ of the the proposed algorithms

As can be seen from Figs. 2 and 4, the proposed algorithms maintain the desired performance for the aforementioned ACI signal. In this subsection, we examine the achieved $P_d$ and $P_f$ of the proposed algorithms for different inband-interference-to-noise ratio (IINR) values which is defined as $IINR \triangleq \frac{\mathrm{E}|\int_{-\infty}^{\infty} f^\star(\tau)(A(t-\tau))d\tau|^2}{\mathrm{E}|\int_{-\infty}^{\infty} f^\star(\tau)(w(t-\tau))d\tau|^2}$ (i.e., for different ACI levels). The IINR is controlled by varying $\kappa$ while setting $\sigma_{\tilde{a}_1} = \sigma_{\tilde{a}_2} = \sigma_s = 1mW$ and $\sigma_w = 4mW$ (i.e., $SNR = -12dB$) in an AWGN channel. For different values of $\kappa$, the equivalent IINR values are obtained by computer simulation and are summarized in **Table I**. As can be seen from Fig. 5, the proposed algorithms maintain the desired $P_f(P_d)$ in the presence of low (moderate) ACI signals. And the performance of the proposed algorithms degrade when the ACI level is high.

We would like to mention here that for any SNR (noise variance), getting the maximum permitted ACI signal level analytically is beyond the scope of this work.

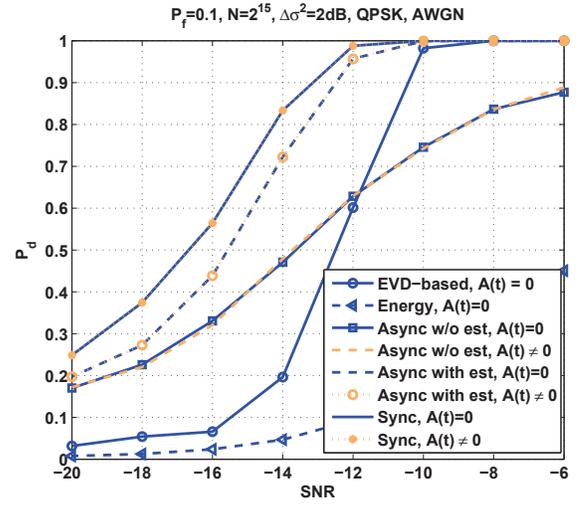

(a)

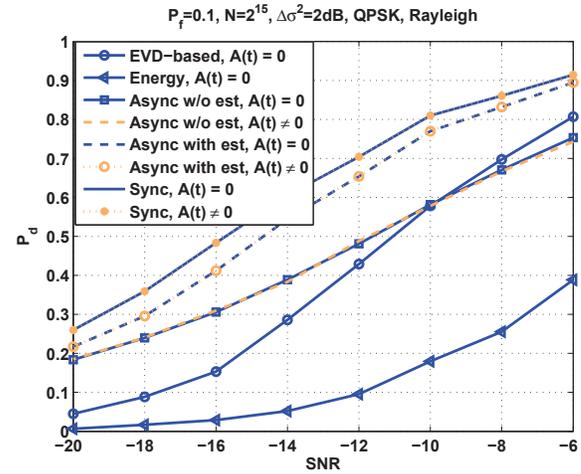

(b)

Fig. 4. Comparison of the proposed detection algorithms, the detection algorithm in [17], EVD-based and energy detection algorithms in (a) AWGN channel. (b) Rayleigh fading channel.

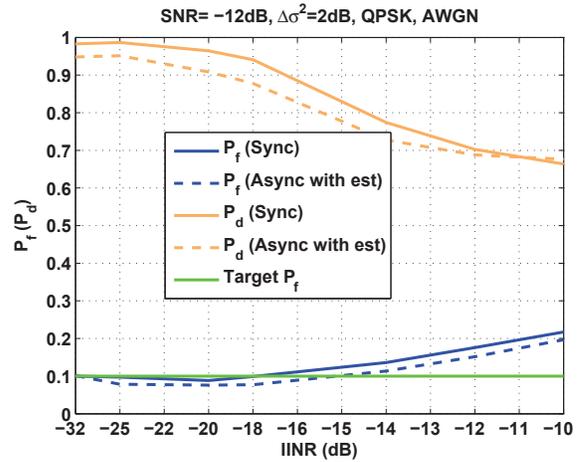

Fig. 5. The effect of ACI on the $P_f$ and $P_d$ of the proposed algorithms.

Table I: $\kappa$ and its equivalent IINR (in dB)

| $\kappa$ | 64 | 128 | 192 | 256 | 320 | 384 | 448 | 512 | 576 | 640 | 704 | 768 |
|---|---|---|---|---|---|---|---|---|---|---|---|---|
| $IINR(dB)$ | -32 | -25 | -22.2 | -20 | -17.9 | -16.4 | -15 | -13.8 | -13.3 | -12 | -11.2 | -10 |

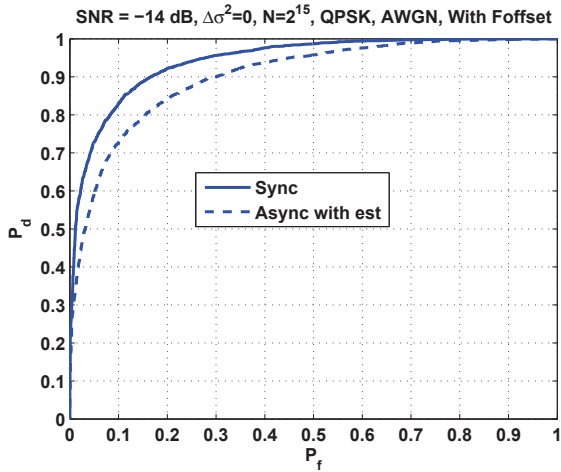

Fig. 6. The performance of the proposed algorithms in the presence of carrier frequency offset.

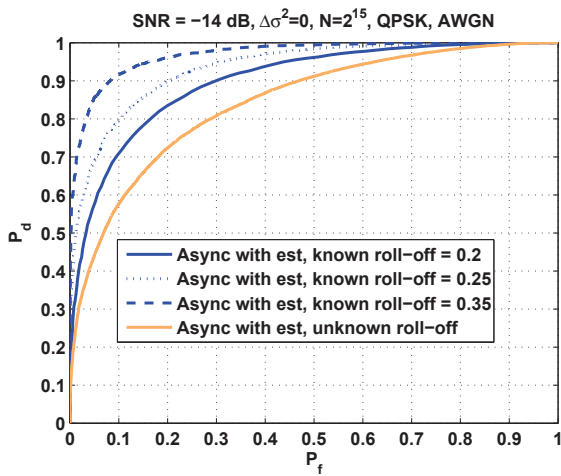

Fig. 7. The performance of the proposed algorithms with known and unknown roll-off factors.

*E. The effect of carrier frequency offset on the performance of the proposed algorithms*

In this simulation, the effect of carrier frequency offset on the detection performances of the proposed detection algorithms are studied. To this end, the carrier frequency offset is modeled as a uniform random variable with interval $\mathcal{U}(-\frac{2*10^{-5}}{T_s}, \frac{2*10^{-5}}{T_s})$. Fig. 6 shows the detection performances of the proposed detection algorithms under this frequency offset and AWGN channel with $A(t) = 0$. From Fig. 3 and Fig. 6, we can observe that carrier frequency offset does not have any impact on the performance of the proposed detection algorithms.

*F. The effect of roll-off factor on the performance of the proposed algorithms*

In this simulation, we examine the effect of roll-off factor on the performance of the proposed algorithms. We also examine the performance of the algorithm of Section III-B (i.e., the detection algorithm for unknown pulse shaping filter). To this end, we consider a DVB-S2 signal which employs a SRRCF with roll-off factor $0.2, 0.25$ or $0.35$ all with equal probability. For this signal, we examine the detection performances of the algorithms of Sections III-A and III-B (i.e., the detection algorithms for known and unknown pulse shaping filters) in AWGN channel with $A(t) = 0$ which is shown in Fig. 7. From this figure, we can observe that better detection performance can be achieved by exploiting the knowledge of the transmitter pulse shaping filter. Furthermore, increasing the roll-off factor increases the detection performance which is expected.

## VI. CONCLUSIONS

This paper proposes novel spectrum sensing algorithms for cognitive radio networks. By assuming known transmitter pulse shaping filter, synchronous and asynchronous receiver scenarios have been considered. For each of these scenarios, the proposed detection algorithm is explained as follows: First, by introducing a combiner vector, the over-sampled signal of total duration equal to $T_s$ is combined linearly. Second, for this combined signal, the SNR maximization and minimization problems are formulated as Rayleigh quotient optimization problems. Third, by employing the solutions of these problems, the ratio of the signal energy corresponding to the maximum and minimum SNRs are proposed as the test statistics. For this test statistics, analytical $P_f$ and $P_d$ expressions are derived for an AWGN channel. The generalization of the proposed algorithms for unknown transmitter pulse shaping filter scenarios has also been discussed. As the $P_f$ expressions do not depend on the true noise variance, the proposed algorithms are robust against noise variance uncertainty. The theoretical expressions are confirmed by computer simulations. Under noise variance uncertainty, simulation results demonstrate that the proposed detection algorithms achieve better detection performance compared to that of the EVD-based and energy detection algorithms in AWGN and Rayleigh fading channels for both synchronous and asynchronous receiver scenarios. Furthermore, simulation results show that the proposed algorithms maintain the desired $P_f(P_d)$ in the presence of low (moderate) ACI signals.

## APPENDIX A
### COMPUTATION OF $\boldsymbol{\Phi}$

The coefficients of $\boldsymbol{\Phi}$ are expressed as

$$\boldsymbol{\Phi}_{(1,1)} = N\mathrm{E}\{\widehat{M}_{a2z}\widehat{M}_{a2z}\} - NM_{a2z}^2$$
$$\boldsymbol{\Phi}_{(1,2)} = \boldsymbol{\Phi}_{(2,1)} = N\mathrm{E}\{\widehat{M}_{a2z}\widehat{M}_{a2e}\} - NM_{a2z}M_{a2e}$$
$$\boldsymbol{\Phi}_{(2,2)} = N\mathrm{E}\{\widehat{M}_{a2e}\widehat{M}_{a2e}\} - NM_{a2e}^2. \quad (34)$$

As can be seen from (14), since the size of $\boldsymbol{\eta}(\boldsymbol{\theta})$ is larger than $L$, $z[n](e[n])$ and $(z[n\pm p], e[n\pm p])$, $\exists p$ are correlated. If $P$ consecutive samples of $z[n](e[n])$ are correlated, by applying multivariate central limit theorem [14], the coefficients of $\boldsymbol{\Phi}$ can be reexpressed as

$$\boldsymbol{\Phi}_{(1,1)} = \mathrm{E}\{|z[n]|^2 \sum_{p=0}^{P} |z[n\pm p]|^2\} - (2P+1)M_{a2z}^2$$

$$\boldsymbol{\Phi}_{(1,2)} = \boldsymbol{\Phi}_{(2,1)} = \mathrm{E}\{|z[n]|^2 \sum_{p=0}^{P} |e[n\pm p]|^2\} - (2P+1)M_{a2z}M_{a2e}$$

$$\boldsymbol{\Phi}_{(2,2)} = \mathrm{E}\{|e[n]|^2 \sum_{p=0}^{P} |e[n\pm p]|^2\} - (2P+1)M_{a2e}^2.$$

Since $z[n]$ and $(z[n\pm p], e[n\pm p])$ are correlated, closed form expression for the coefficients of $\boldsymbol{\Phi}$ is not possible. In the following, we provide a numerical method to compute the expectation $\mathrm{E}\{|z[n]|^2|e[n+p]|^2\}$ (all the other expectation terms can be computed similar to this expectation).

By defining $|z[n]|^2 = z_r[n]^2 + z_i[n]^2$ and $|e[n]|^2 = e_r[n]^2 + e_i[n]^2$, where $(.)_r$ and $(.)_i$ denote real and imaginary, respectively, it can be shown that

$$\mathrm{E}\{|z[n]|^2|e[n+p]|^2\} =$$
$$\mathrm{E}\{(z_r[n]e_r[n+p])^2 + (z_i[n]e_i[n+p])^2\}$$
$$+ \mathrm{E}\{z_r[n]^2 e_i[n+p]^2 + z_i[n]^2 e_r[n+p]^2\}$$
$$= 2\sigma_z^2 \sigma_e^2 (\mathrm{E}\{(\bar{z}[n]\bar{e}[n+p])^2\} + 1) \quad (35)$$

where $\bar{e}[n](\bar{z}[n]) \sim \mathcal{N}(0,1)$.

As we can see, by applying zero padding, it is possible to express $\bar{z}[n]$ and $\bar{e}[n+p]$ as a fully correlated samples of appropriate size. To this end, we derive an expression for the expectation term $\mathrm{E}\{(\tau\eta)^2\}$, where

$$\tau = \sum_{i=1}^{J} c_i \tilde{w}_i, \quad \eta = \sum_{i=1}^{J} d_i \tilde{w}_i$$

$J$ is a positive integer, $\{c_i, d_i\}_{i=1}^{J}$ are arbitrary coefficients and $\{\tilde{w}_i\}_{i=1}^{J}$ are i.i.d zero mean Gaussian random variables all with unit variance.

$$\mathrm{E}\{(\tau\eta)^2\} = \mathrm{E}\{(\sum_{i=1}^{J} c_i d_i \tilde{w}_i^2 + \sum_{i=1}^{J} \sum_{j=1, j\neq i}^{J} c_i d_j \tilde{w}_i \tilde{w}_j)^2\}$$
$$= \mathrm{E}\{(\sum_{i=1}^{J} c_i d_i \tilde{w}_i^2)^2 + (\sum_{i=1}^{J} \sum_{j=i, j\neq i}^{J} c_i d_j \tilde{w}_i \tilde{w}_j)^2\}$$
$$= \sum_{i=1}^{J} \sum_{j=1}^{J} c_i^2 d_j^2 + 2c_i c_j d_i d_j \quad (36)$$

where the second equality employs $\mathrm{E}\{\tilde{w}_i \tilde{w}_j = 0, \forall i \neq j\}$ and the third equality employs the definition of a fourth moment of a Gaussian random variable [14].

## APPENDIX B
### COMPUTATION OF $\{\boldsymbol{\Sigma}_k\}_{k=1}^{L}$

In this appendix, we provide numerical methods to compute the covariance matrix of $\boldsymbol{\Sigma}_1$. By extending *Theorem 1* to multivariate Gaussian random variables, it can be shown that [12]

$$\mathbf{t}_1 = [T_1, T_2, \cdots, T_L]^T = \quad (37)$$
$$\sqrt{N}[(\widehat{\widetilde{T}}_1 - \widetilde{T}_1), \cdots, (\widehat{\widetilde{T}}_L - \widetilde{T}_L)]^T \sim \mathcal{N}(\mathbf{0}, \boldsymbol{\Sigma}_1)$$

where $\boldsymbol{\Sigma}_1 = \widetilde{\mathbf{V}}\widetilde{\boldsymbol{\Phi}}\widetilde{\mathbf{V}}^T$,

$$\widetilde{\mathbf{V}}_{(i,:)} = \left[\frac{\partial \widehat{\widetilde{T}}_i}{\partial \widehat{M}_{a2z1}}, \frac{\partial \widehat{\widetilde{T}}_i}{\partial \widehat{M}_{a2e1}}, \frac{\partial \widehat{\widetilde{T}}_i}{\partial \widehat{M}_{a2z2}}, \frac{\partial \widehat{\widetilde{T}}_i}{\partial \widehat{M}_{a2e2}}, \cdots, \right.$$
$$\left. \frac{\partial \widehat{\widetilde{T}}_i}{\partial \widehat{M}_{a2zL}}, \frac{\partial \widehat{\widetilde{T}}_i}{\partial \widehat{M}_{a2eL}}\right]_{\{\widehat{M}_{a2zj}=M_{a2zj}, \widehat{M}_{a2ej}=M_{a2ej}\}_{j=1}^{L}},$$

and $\widetilde{\boldsymbol{\Phi}}$ is the asymptotic covariance matrix of a multivariate random variable $\sqrt{N}[\widehat{M}_{a2z1} - M_{a2z1}, \widehat{M}_{a2e1} - M_{a2e1}, \widehat{M}_{a2z2} - M_{a2z2}, \widehat{M}_{a2e2} - M_{a2e2}, \cdots, \widehat{M}_{a2zL} - M_{a2zL}, \widehat{M}_{a2eL} - M_{a2zL}]^T \sim \mathcal{N}(\mathbf{0}, \widetilde{\boldsymbol{\Phi}})$ with $\widehat{M}_{a2zi}(\widehat{M}_{a2ei})$ as the estimated absolute second moment corresponding to $\widehat{\widetilde{T}}_i$.

The coefficients of $\widehat{\mathbf{V}}$ can be computed analytically like that of (18) and the coefficients of $\widetilde{\boldsymbol{\Phi}}$ can be computed numerically like that of (34). The details are omitted for conciseness. Note that the covariance matrices $\{\boldsymbol{\Sigma}_k\}_{k=2}^{L}$ can be computed numerically like that of $\boldsymbol{\Sigma}_1$.

## APPENDIX C
### STRATEGY FOR DETERMINING $\mathbf{A}$, $\mathbf{B}$ AND $\boldsymbol{\alpha}$

When $f^{\star}(t) = g(t)$ is SRRCF, it can be easily seen that

$$\int_{-\infty}^{\infty} f^{\star}(\tau)f(t+\tau)d\tau = \tilde{f}(t)$$

where $\tilde{f}(t) = h(t)$ is a raised cosine filter (RCF). It follows from (5) that $\mathbf{B}_{(i+1,j+1)} = \tilde{f}(t_i - t_j), \forall i, j$. When the roll-off factor = 0.2 and L = 8, we will get (38).

In general, we are not aware of any analytical approach to compute all the entries of $\mathbf{A}$. However, some of the entries of this matrix can be computed in closed form as follows: From the properties of RCF, we have $h(t) = 0$, for $t = kT_s, k = 1, 2, \cdots$. By exploiting this property and setting $t_0 = \frac{T_s}{2}$ (i.e., as on page 5), it can be shown that $\mathbf{A}_{(t_0+1,:)} = \mathbf{B}_{(t_0+1,:)}$.

From this explanation, we can understand that most of the entries of $\mathbf{A}$ are computed numerically. To compute $\mathbf{A}$ numerically, we need to set the maximum $k'$ (i.e., $K'$). From numerical computation, we have observed that each entry of $\mathbf{A}$ is almost constant when $K' \geq 2^{15}$. For this reason, we employ $K' = 2^{15}$ in our simulation. With this $K'$, the entries of $\mathbf{A}$ are given by (39).

By applying the above $\mathbf{A}$ and $\mathbf{B}$, the optimal $\boldsymbol{\alpha}_{min}$ of (6) and $\boldsymbol{\alpha}_{max}$ of (7) can be expressed as in (40).

In this example, (10) is computed by replacing very small Eigenvalues (i.e., $< 10^{-3}$) of $\mathbf{B}$ with 0.

$$\mathbf{B} = \begin{bmatrix} 1.0000 & 0.9739 & 0.8982 & 0.7801 & 0.6307 & 0.4637 & 0.2938 & 0.1353 \\ 0.9739 & 1 & 0.9739 & 0.8982 & 0.7801 & 0.6307 & 0.4637 & 0.2938 \\ 0.8982 & 0.9739 & 1 & 0.9739 & 0.8982 & 0.7801 & 0.6307 & 0.4637 \\ 0.7801 & 0.8982 & 0.9739 & 1 & 0.9739 & 0.8982 & 0.7801 & 0.6307 \\ 0.6307 & 0.7801 & 0.8982 & 0.9739 & 1 & 0.9739 & 0.8982 & 0.7801 \\ 0.4637 & 0.6307 & 0.7801 & 0.8982 & 0.9739 & 1 & 0.9739 & 0.8982 \\ 0.2938 & 0.4637 & 0.6307 & 0.7801 & 0.8982 & 0.9739 & 1 & 0.9739 \\ 0.1353 & 0.2938 & 0.4637 & 0.6307 & 0.7801 & 0.8982 & 0.9739 & 1 \end{bmatrix}. \quad (38)$$

$$\mathbf{A} = \begin{bmatrix} 0.9000 & 0.8816 & 0.8276 & 0.7420 & 0.6307 & 0.5016 & 0.3635 & 0.2258 \\ 0.8816 & 0.9146 & 0.9086 & 0.8629 & 0.7801 & 0.6658 & 0.5284 & 0.3780 \\ 0.8276 & 0.9086 & 0.9500 & 0.9469 & 0.8982 & 0.8071 & 0.6804 & 0.5284 \\ 0.7420 & 0.8629 & 0.9469 & 0.9854 & 0.9739 & 0.9128 & 0.8071 & 0.6658 \\ 0.6307 & 0.7801 & 0.8982 & 0.9739 & 1.0000 & 0.9739 & 0.8982 & 0.7801 \\ 0.5016 & 0.6658 & 0.8071 & 0.9128 & 0.9739 & 0.9854 & 0.9469 & 0.8629 \\ 0.3635 & 0.5284 & 0.6804 & 0.8071 & 0.8982 & 0.9469 & 0.9500 & 0.9086 \\ 0.2258 & 0.3780 & 0.5284 & 0.6658 & 0.7801 & 0.8629 & 0.9086 & 0.9146 \end{bmatrix}. \quad (39)$$

$$\begin{aligned} \boldsymbol{\alpha}_{min} &= \begin{bmatrix} -8.8565 & 5.2981 & 8.3685 & 3.8283 & -3.4689 & -8.1746 & -5.4128 & 8.3317 \end{bmatrix} \\ \boldsymbol{\alpha}_{max} &= \begin{bmatrix} -0.0586 & -0.0033 & 0.1166 & 0.2486 & 0.3346 & 0.3213 & 0.1697 & -0.1375 \end{bmatrix}. \end{aligned} \quad (40)$$